\documentstyle[12pt,theorem,amssymb,latexsym]{article}
\advance \topmargin by -\headheight
\advance \topmargin by -\headsep
\evensidemargin \oddsidemargin
\author{Jean-Paul Allouche \\
{\small CNRS, LRI} \\
{\small Universit\'e Paris-Sud, B\^at. 490} \\
{\small F-91405 Orsay Cedex (France)} \\
{\small\tt allouche@lri.fr} \\
\and
Michael Baake \\
{\small Institut f\"ur Mathematik und Informatik} \\
{\small Universit\"at Greifswald, Jahnstra{\ss}e 15 a} \\
{\small D-17487 Greifswald (Germany)} \\
{\small\tt mbaake@uni-greifswald.de} \\
\and
Julien Cassaigne \\
{\small CNRS, IML} \\
{\small Case 907, 163 Avenue de Luminy} \\
{\small F-13288 Marseille Cedex 9 (France)} \\
{\small\tt cassaigne@iml.univ-mrs.fr} \\
\and
David Damanik \\
{\small Department of Mathematics 253-37} \\
{\small California Institute of Technology} \\
{\small Pasadena, CA 91125 (USA)} \\
{\small\tt damanik@its.caltech.edu} \\
}

\title{Palindrome complexity}

\date{ }

\def \proof{\bigbreak\noindent{\it Proof.\ \ }}

\def \endpf{{\ \ $\Box$}}

\newtheorem{theorem}{Theorem}
\newtheorem{lemma}{Lemma}

\theorembodyfont{\rm}

\newtheorem{remark}{Remark}

\newtheorem{definition}{Definition}

\begin{document}

\maketitle

\begin{center} 
Dedicated to Jean Berstel for his 60th birthday with our very best wishes.
\end{center}

\begin{abstract}
We study the {\em palindrome complexity\,} of infinite sequences on finite 
alphabets, i.e., the number of palindromic factors (blocks) of given 
length occurring in a given sequence. We survey the known results and obtain 
new results for some sequences, in particular for Rote sequences and for 
fixed points of primitive morphisms of constant length belonging to 
``class P'' of Hof-Knill-Simon.
We also give an upper bound for the palindrome complexity of a sequence in
terms of its (block-)complexity.
\end{abstract}

\section{Introduction} The (block-)complexity function of an infinite 
sequence on a finite alphabet is the number of factors (blocks) of given 
length occurring in this sequence. This notion was introduced in 1975 by 
Ehrenfeucht, Lee, and Rozenberg \cite{ELR}. The complexity function of a
sequence measures in some sense how ``complicated'' the sequence is: the 
reader is referred to the surveys \cite{Allcomp,Fer2,FerKas}.
Due {\em inter alia\,} to its applications to physics \cite{HKS}, another
interesting ``complexity'' for an infinite sequence on a finite alphabet 
is its {\em palindrome complexity}, i.e., the number of palindromic factors 
(blocks) of given length occurring in the sequence. Combinatorial results 
on palindrome complexity for some sequences or classes of sequences were
proved in \cite{Dro,Allschro,DroPir,Baa,Dam1,DamZar}. We survey known 
results, and give new ones, proving for example that {\em the palindrome 
complexity of Rote sequences is constant and equal to $2$}. We prove that
{\em the palindrome complexity function of a sequence that is a fixed point 
of a primitive morphism of length $d$ belonging to ``class P'' of 
Hof-Knill-Simon and satisfying some technical conditions is $d$-automatic}. 
Finally, we give {\em an upper bound for the palindrome complexity in terms 
of the (usual) complexity}. 

\section{Definitions and notations}

\subsection{Generalities}

We will use the notation ${\Bbb N} = \{0, 1, 2, \ldots\}$.
We recall the following classical definitions in combinatorics of words.
The sequences we consider in this paper are defined on a finite 
{\sl alphabet\,} (i.e., on a finite set) ${\cal A}$. Elements of 
${\cal A}$ are called {\sl letters}. The set ${\cal A}^*$ is defined as
the set of {\sl words\,} on ${\cal A}$, i.e., the set of (possibly empty) 
strings of symbols of ${\cal A}$, equipped with the {\it concatenation}.
In other words, ${\cal A}^*$ is the free monoid for the concatenation 
generated by ${\cal A}$. The {\sl length\,} of a word $w$, denoted by $|w|$, 
is recursively defined by: the empty word has length $0$, and for any word 
$w$ and any letter $a$, $|wa| = |w| + 1$. A word is called a {\sl factor\,} 
of another word or of an infinite sequence if it occurs ``without hole'' in
this word or sequence (factors are also called {\sl subwords}, while the 
term {\sl substrings\,} stands for the case where there are holes: $001$ is 
a factor of $11001110$ but only a substring of $0101010$). 

If ${\cal A}$ and ${\cal B}$ are two alphabets, homomorphisms for the
concatenation from ${\cal A}^*$ to ${\cal B}^*$ are called {\sl morphisms}.
A {\sl morphism\,} on ${\cal A}$ (also called {\sl substitution\,} or 
{\sl inflation rule}) is a morphism from ${\cal A}^*$ into itself. 
A morphism is defined by its values on the letters. A {\sl constant length 
morphism\,} (or {\sl uniform morphism}) is a morphism such that the images of 
all letters have the same length. If this length is equal to $d$, the morphism 
is also called a {\it morphism of length $d$}, or a {\it $d$-morphism}. 
A sequence on ${\cal A}$ is called {\sl $d$-automatic} if it is the pointwise 
image of a fixed point of a morphism of length $d$ on an alphabet ${\cal B}$ 
(pointwise image means, of course, image under a morphism of length $1$ from 
${\cal B}$ to ${\cal A}$).
A morphism $\varphi$ on ${\cal A}$ is called {\em primitive\,} if there exists 
an integer $k \geq 1$ such that the image of each letter by $\varphi^k$
contains at least one occurrence of each letter of ${\cal A}$.
A morphism is called {\em non-erasing\,} if the image of every letter is
different from the empty word.
 
A (non-erasing) morphism on ${\cal A}$ can be extended to infinite sequences 
with values in ${\cal A}$ ``by continuity''. (The set of sequences 
${\cal A}^{\Bbb N}$ is equipped with the topology of simple convergence, 
i.e., the product topology where each copy of ${\cal A}$ is equipped with the 
discrete topology.) 
An infinite sequence can thus be a {\em fixed point\,} of a morphism. 

An infinite sequence on a finite alphabet is called {\em recurrent\,} if
each word that occurs in the sequence occurs infinitely often. The sequence
is called {\em uniformly recurrent\,} or {\em minimal\,} if each word that
occurs in the sequence occurs infinitely often and the distance between two
consecutive occurrences is bounded (some authors use the term
{\em almost-periodic\,} for such sequences, while some authors call them
{\em repetitive}). It is easy to prove that a sequence that is a fixed point 
of a primitive morphism is uniformly recurrent.

\subsection{Periodicity}

In this section we recall the definition of periodic sequences or words,
and we give two theorems and a lemma that will prove useful.

\begin{definition}\label{period}

\mbox{ }

\begin{itemize}

\item An infinite sequence $u = u_0 u_1 \ldots$ is called {\em periodic\,} if 
there exists an integer $T \geq 1$ (called a {\em period} of the sequence) 
such that for each $n \geq 0$ we have $u_{n+T} = u_n$. It is called 
{\em ultimately periodic\,} if there exists an integer $\ell$ such that the 
sequence $v$ defined by $v_n = u_{\ell + n}$ is periodic.

\item Let $w$ be a finite word. Any integer $p \geq 1$ such that $w$ is 
a prefix of an infinite sequence of period $p$ is called {\em a period\,} 
of $w$. {\em The period\,} of $w$ is the smallest such integer.
(For example a period of the word $01101$ is  $5$. Its period is $3$.)
\end{itemize}

\end{definition}

Theorems \ref{finewilf} and \ref{lynsch} below will prove useful. They are 
due respectively to Fine and Wilf \cite{FinWil} and to Lyndon and 
Sch\"utzenberger~\cite{LynSch}.

\begin{theorem}[Fine-Wilf]\label{finewilf}
Let $u = (u_n)_{n \geq 0}$ and $v = (v_n)_{n \geq 0}$ be two periodic 
sequences with respective periods $T$ and $T'$. If $u_n = v_n$ for more than
$T + T' - \gcd(T,T')$ consecutive values of $n$, then the sequences $u$ 
and $v$ are equal. The value $T + T' - \gcd(T,T')$ is sharp.
\end{theorem}

\begin{remark}
In the literature an easy corollary of this result is often called the
theorem of Fine and Wilf, namely that {\em if a finite word $w$ has
two periods $T$ and $T'$ such that $|w| \geq T + T' - \gcd(T,T')$, then 
$\gcd(T,T')$ is a period of $w$. The value $T + T' - \gcd(T,T')$ is sharp}.
\end{remark}

\begin{theorem}[Lyndon-Sch\"utzenberger]\label{lynsch}
Let ${\cal A}$ be an alphabet. Let $x,y,z \in {\cal A}^*$, with
$x$ and $z$ non-empty.  Then $xy = yz$ if and only if
there exist $u,v \in {\cal A}^*$, and an integer $e \geq 0$ such that
$x = uv$, $z = vu$, and
$y = (uv)^eu = u(vu)^e$.
\end{theorem}

\begin{lemma}\label{periods}

\mbox{ }

\begin{itemize}

\item If the period $T$ of a word $w$ satisfies $T \leq |w|/2$, then all 
the periods of $w$ that are $\leq |w|/2$ are divisible by $T$.

\item Let $z$ be a word and let $w$ be a factor of $z$. If $z$ has period
$T$, if $w$ has period $T'$, and if $T + T' \leq |w|$, then $T'$ is a 
period of $z$.

\item Let $z$ be a word and let $w$ and $w'$ be two factors of $z$. 
Suppose that $T$, the period of $w$, $T'$, the period of $w'$, 
and $\Theta$, the period of $z$, satisfy 
$\Theta + T \leq |w|$ and $\Theta + T' \leq |w'|$.
Then $T = T'$.

\end{itemize}

\end{lemma}

\proof

\begin{itemize}

\item If $T$ is the period of $w$, if $U$ is another period, then $w$ is a 
prefix of a $T$-periodic sequence and of a (possibly distinct) $U$-periodic 
sequence. If both periods are $\leq |w|/2$, then $T+U \leq |w|$. Hence
the two infinite periodic sequences, coinciding on a prefix of length 
$\geq T+U$, must be equal (from Theorem~\ref{finewilf}).
The periodic sequence thus obtained admits $T$ as smallest period (namely
$T$ is the smallest period of the word $w$). Hence $T$ divides $U$.
(Actually it can also be noted that if $T > |w|/2$, then our statement is 
empty, hence holds.)

\item The word $z$ is a prefix of a $T$-periodic sequence $u$, and the word
$w$ is a prefix of a $T'$-periodic sequence $u'$. But the word $w$ is a factor 
of $z$, hence a prefix of a sequence of period $T$, say $\hat{u}$, obtained 
from the sequence $u$ by erasing some prefix. Since $T+T' \leq |w|$ we 
have, from Theorem~\ref{finewilf}, that $u' = \hat{u}$. 
Hence $T'$ is a period of $u$, since $u$ can be obtained by erasing a prefix 
of $\hat{u}$. Hence $T'$ is a period of $z$.

\item From the second item above, $T$ and $T'$ are periods of $z$. Hence
clearly $T'$ is a period of $w$ and $T$ is a period of $w'$. 
Since $T$ is the minimal period of $w$ and $T'$ the minimal period of $w'$,
we have $T\leq T'$ and $T' \leq T$.
Hence $T'=T$. \endpf
\end{itemize}

\subsection{Palindromes and complexity}

\begin{definition}
If $w = w_1 w_2 \ldots w_j$ is a word on the alphabet ${\cal A}$,
we denote by $\widetilde{w}$ the word obtained by reading $w$ backwards, 
i.e., $\widetilde{w} = w_j w_{j-1} \ldots w_2 w_1$. A {\em palindrome\,} 
is a word $w$ such that $w = \widetilde{w}$. (For example the words ``level''
and ``deed'' are palindromes in English.)
\end{definition}

\begin{definition}
Let $u := u_0 u_1 u_2 \ldots$ be a sequence on the finite alphabet ${\cal A}$.
We denote by $\mbox{\rm fac}_u(n)$ the number of words of length $n$ that 
are factors of the sequence $u$. We denote by $\mbox{\rm pal}_u(n)$ the number 
of palindromes of length $n$ that are factors of the sequence $u$. 
\end{definition}

\begin{remark}

The notations used in other papers might differ. In the literature
$p_u(n)$ sometimes stands for the block-complexity and sometimes for
the palindrome complexity. We hope that our terminology is unambiguous.

\end{remark}

\subsection{Sturmian sequences}

We end this section by recalling the definition of Sturmian sequences.
These sequences can be obtained by playing billiard on squares,
starting with an irrational slope. They can also be defined by their
complexity. The reader is referred to \cite{MorHed,CovHed,LunPle,Ber1,Ber2}.

\begin{definition}
A {\em Sturmian sequence\,} is a sequence $u = (u_n)_{n \geq 0}$ whose
(block-)complexity ${\rm fac}_u(k)$ satisfies: $\forall k \geq 1$,
${\rm fac}_u(k)= k+1$.
\end{definition}

\section{Motivation in physics}

Given a uniformly recurrent sequence $u$, one may consider the associated
{\it LI-class}~$\Omega$ (also called {\it hull} or {\it induced subshift}) 
which consists of all two-sided infinite sequences that have the same finite 
factors as $u$. If $u$ is Sturmian, if $u$ is generated by a primitive 
morphism, or if $u$ is derived from a standard cut and project scheme, this 
gives widely used models of {\it one-dimensional quasicrystals}. That is, a 
one-dimensional quasicrystal is modelled by a suitable family of two-sided 
sequences which are locally indistinguishable since they have the same finite 
factors.

To such a structure, one may associate a family of discrete one-dimensional
Schr\"odinger operators $(H_\omega)_{\omega \in \Omega}$ as follows: choose 
an injective function $f : {\cal A} \rightarrow {\Bbb R}$ and define, for 
every $\omega \in \Omega$, the operator $H_\omega$ in $\ell^2({\Bbb Z})$ by
$$
(H_\omega \phi)(n) = \phi(n+1) + \phi(n-1) + f(\omega_n) \phi(n).
$$
The spectral properties of $H_\omega$ determine the ``conductivity 
properties'' of the given structure. Roughly speaking, if the spectrum is 
absolutely continuous, then the structure behaves like a conductor, while 
in the case of pure point spectrum, it behaves like an insulator. 
The intermediate spectral type -- singular continuous spectrum -- is 
generally expected to give rise to intermediate transport properties, but 
no one currently understands this correspondence very well.

For classical (periodic or atomic) structures, singular continuous spectra 
do not occur. However, for one-dimensional quasicrystals, this spectral type 
appears to be rather typical, and there has been a lot of recent research 
activity focussing on results of this kind. One important contribution in 
this direction is the paper of Hof, Knill, and Simon \cite{HKS}, 
where a sufficient criterion for purely singular continuity of individual 
such operators is derived in terms of a strong palindromicity property of 
the underlying sequence.

 From a combinatorial point of view, it is interesting to note that the 
criterion of Hof, Knill, and Simon deduces singular continuous spectrum from 
explicit combinatorial (and dynamical) properties of $u$. Suppose that the 
factors of $u$ occur with well-defined, positive frequencies ($\Omega$ 
is strictly ergodic). This is the case for a large class of sequences, 
including Sturmian sequences and sequences generated by primitive morphisms, 
but also for sequences derived from the standard or generalized cut and 
project schemes \cite{Sch}. Then the following holds \cite{HKS}.

\begin{theorem}[Hof-Knill-Simon]
Let $u$ be a sequence on a finite alphabet that is not ultimately 
periodic (hence ${\rm fac}_u(k) \ge k+1$ for every $k \geq 1$), 
and such that ${\rm pal}_u$ does not ultimately vanish, i.e., 
$\limsup_{k \rightarrow \infty} {\rm pal}_u(k) > 0$. We suppose 
that the subshift $\Omega$ induced by $u$ is strictly ergodic. 
Then, for uncountably many $\omega \in \Omega$, the operator 
$H_\omega$ has purely singular continuous spectrum.
\end{theorem}

This theorem applies to a large class of sequences generated by primitive 
morphisms, to all Sturmian sequences, and, more generally, to all sequences 
defined by circle maps. Furthermore, it also applies to sequences derived from 
a standard cut and project scheme with inversion-symmetric window (see
\cite{Baa} for details).
Note that by uniform recurrence of $u$, it suffices to assume that $u$ is not
periodic. We remark that there is a similar, purely combinatorial, sufficient
condition for purely singular continuous spectrum in terms of {\em powers\,}
occurring in $u$, see \cite{Dam2} for a survey.

\section{Survey of known results}

\subsection{Rudin-Shapiro and paperfolding sequences}

Paperfolding sequences are binary sequences obtained by repeatedly folding
a strip of paper (see \cite{DMFP} for example). We recall the definition.

\begin{definition}

A sequence $(u_n)_{n \geq 1}$ with values in $\{0,1\}$
is called a {\em paperfolding sequence\,} if there exists a sequence
$i_0, i_1, i_2, \ldots$, with $i_k \in \{0, 1\}$ (called the sequence
of {\em unfolding instructions}) such that
$$
\forall m \geq 0, \ \forall j \geq 0, \ \
u_{2^m(2j+1)} \equiv j + i_m \bmod 2.
$$
(Note that every integer $n \geq 1$ can be uniquely written as
$n = 2^m(2j+1)$, with $m, j \geq 0$.)

\end{definition}

\bigskip

Generalized Rudin-Shapiro sequences (in the sense of \cite{MenTen}) are
obtained by ``integrating modulo $2$'' paperfolding sequences.
More precisely

\begin{definition}\label{gRS}

A sequence $(v_n)_{n \geq 0}$ is called a {\em generalized Rudin-Shapiro 
sequence\,} if $v_0 = 0$ and if there exists a paperfolding sequence 
$(u_n)_{n \geq 1}$ such that, for every $n \geq 0$, there holds
$v_n \equiv \sum_{k=1}^n u_k \bmod 2$. (The ``classical'' Rudin-Shapiro 
sequence corresponds to the paperfolding sequence with unfolding 
instructions $0, 0, 1, 0, 1, 0, 1, 0, 1, \ldots$)

\end{definition}

\begin{remark}
There exist other generalizations of the classical Rudin-Shapiro sequence 
(see, e.g., \cite{Que,AllLia}) but we restrict here to the ones given in 
Definition~\ref{gRS} above.
\end{remark}

The following theorem was proved in \cite{Allschro}. It was proved again in a 
more efficient and general way by Baake in \cite{Baa}.

\begin{theorem}[Allouche]\label{RS}

\mbox{ }

The palindrome complexity of any paperfolding sequence $u=(u_n)_{n \geq 1}$ 
satisfies ${\rm pal}_u(k) = 0$ for any $k \geq 14$.

The palindrome complexity of any generalized Rudin-Shapiro sequence 
$v=(v_n)_{n \geq 0}$ satisfies ${\rm pal}_v(k) = 0$ for any $k \geq 15$.
\end{theorem}

\subsection{The period-doubling sequence}

\begin{definition}
The {\em period-doubling sequence\,} is defined as the infinite fixed point
of the morphism $0 \to 01$, $1 \to 00$.
\end{definition}

The following result was proved in \cite{Dam1}.

\begin{theorem}[Damanik]
The palindrome complexity of the period-doubling sequence $u$ satisfies
$$
\left\{\begin{array}{ll}
\forall k \ {\rm even}, \ k \geq 4 & \Rightarrow {\rm pal}_u(k) = 0, \\
\forall k \ {\rm odd}, \ k \geq 5 & \Rightarrow {\rm pal}_u(k) 
= {\rm pal}_u(2k-1) = {\rm pal}_u(2k+1),
\end{array}
\right.
$$
and the first few values of \,$\mbox{\rm pal}_u(k)$ are given by
$$
\vbox{\offinterlineskip
\hrule
\halign{
\vrule# & \strut $\quad \hfil# \quad$&
\vrule# & \strut $\quad \hfil# \quad$&
\vrule# & \strut $\quad \hfil# \quad$&
\vrule# & \strut $\quad \hfil# \quad$&
\vrule# & \strut $\quad \hfil# \quad$&
\vrule# & \strut $\quad \hfil# \quad$&
\vrule# & \strut $\quad \hfil# \quad$&
\vrule# & \strut $\quad \hfil# \quad$&
\vrule# & \strut $\quad \hfil# \quad$&\vrule#  \cr
&k\hfil&&1\hfil&&2\hfil&&3\hfil&&4\hfil&&5\hfil&&6\hfil&&7\hfil&\cr
\noalign{\hrule}
&\mbox{\rm pal}_u(k)&&2&&1&&3&&0&&4&&0&&3&\cr
}
\hrule}
$$
In particular, the function $k \to {\rm pal}_u(k)$ takes its values in the set 
$\{0, 1, 2, 3, 4\}$. Furthermore, $\limsup_{k \to \infty} {\rm pal}_u(k) = 4$.
(Actually\, $0$, $3$, and $4$ are the only values that are taken infinitely 
often.)
\end{theorem} 

\subsection{Sturmian sequences}

The following nice characterization of Sturmian sequences in terms of
palindrome complexity was given in \cite{DroPir} (see also \cite{Dro}
for the Fibonacci sequence, which can be defined as the fixed point of
the morphism $0 \to 01$, $1 \to 0$).

\begin{theorem}[Droubay-Pirillo]\label{dropir}

A sequence $u = (u_n)_{n \geq 0}$ is Sturmian if and only if its 
palindrome complexity satisfies: $\forall k$ odd, ${\rm pal}_u(k) = 2$ 
and\, $\forall k \geq 2$ even, ${\rm pal}_u(k) = 1$.

\end{theorem}

\begin{remark}
A generalization of this theorem to the two-dimensional case has been obtained
recently by Berth\'e and Vuillon \cite{BerVui}. For a study of palindromes
in {\sl episturmian sequences}, see \cite{DroJusPir,JusPir}.
\end{remark}

\subsection{Fixed points of primitive morphisms}

Primitive morphisms are often considered because they have dynamical or
combinatorial properties that other morphisms may not have. The following 
result was proved in \cite{DamZar}.

\begin{theorem}[Damanik-Zare]\label{damzar}
The palindrome complexity ${\rm pal}_u(k)$ of a fixed point 
$u = (u_n)_{n \geq 0}$ of a primitive morphism is bounded (hence takes 
only finitely many values).
\end{theorem}

\begin{remark} As a consequence of Theorem~\ref{cassaigne} below, we will
have that the conclusion of Theorem~\ref{damzar} also holds for uniform 
(not necessarily primitive) morphisms.
\end{remark}

\section{Rote sequences}

Sequences of (block-)complexity $2k$ were studied in \cite{Rot}.
In that paper, Rote proved in particular \cite[Theorem 3]{Rot} that
{\sl an infinite $0,1$-sequence $w = (w_n)_{n \geq 0}$ is a 
complementation-symmetric sequence with block-complexity $2k$ if and only
if its first difference (modulo $2$) sequence $\beta = (\beta_n)_{n \geq 0}$, 
is Sturmian, where, for each $n \geq 0$, $\beta_n := w_{n+1} - w_n \bmod 2$}. 
Recall that a complementation-symmetric sequence on a two-letter
alphabet, say ${\cal A} = \{a, b\}$, is a sequence such that for any block 
occurring in it, the block obtained by changing $a$'s into $b$'s and $b$'s 
into $a$'s is also a factor. Using Theorem~\ref{dropir} we can compute the 
palindrome complexity of these sequences.

\begin{theorem}\label{roteseq}
Let $w = (w_n)_{n \geq 0}$ be a complementation-symmetric sequence with
complexity ${\rm fac}_w(k)=2k$ for all $k \geq 1$. Then its palindrome
complexity satisfies ${\rm pal}_w(k) = 2$ for all $k \geq 1$.
\end{theorem}

\proof

Since ${\rm fac}_w(1) = 2$, we see that the sequence $w$ is a binary sequence.
Without loss of generality we may assume that it is a $0,1$-sequence.
Define $\beta = (\beta_n)_{n \geq 0}$ by: $\forall n \geq 0$ 
$\beta_n := w_{n+1} - w_n \bmod 2$. Then we know that the sequence
$\beta$ is Sturmian.

Define for each $k \geq 1$ the maps $\Phi_k$ and $\Psi_k^s$ (where $s = 0,1$)
on words of length $k$ on $\{0,1\}$ by
$$
\begin{array}{llll}
\Phi_k(a_1 a_2 \ldots a_k) &:=& b_1 b_2 \ldots b_{k-1}, \ &\mbox{\rm where} \
b_i : = a_{i+1} - a_i \bmod 2, \\
\Psi_k^s(b_1 b_2 \ldots b_k) &:=& s c_1 c_2 \ldots c_k, &\mbox{\rm where} \
s = 0, 1 \ \mbox{\rm and} \ c_j := s + \displaystyle\sum_{i=1}^j b_i \bmod 2.
\end{array}
$$

It is straightforward if $k$ is even that $\Phi_k$ sends palindromes of length
$k$ to palindromes of length $k-1$ whose central letter is $0$, and that any
palindrome $\pi$ of length $k-1$ whose central letter is $0$ is the image under
$\Phi_k$ of exactly two palindromes of length $k$, namely $\Psi_{k-1}^0(\pi)$
and $\Psi_{k-1}^1(\pi)$.

It is also straightforward if $k$ is odd that $\Phi_k$ sends palindromes of 
length $k$ to palindromes of length $k-1$, and that any palindrome $\pi$ of 
length $k-1$ is the image under $\Phi_k$ of exactly two palindromes of length 
$k$, namely $\Psi_{k-1}^0(\pi)$ and $\Psi_{k-1}^1(\pi)$.

Now from this property of the map $\Phi_k$ and from the relation between the 
sequences $w$ and $\beta$, we see that the number of palindromes of length $k$ 
occurring in the sequence $w$ is equal to twice the number of palindromes of 
length $k-1$ whose central letter is $0$ occurring in $\beta$ if $k$ is even, 
and to twice the number of palindromes of length $k-1$ occurring in $\beta$ 
if $k$ is odd. 

To conclude the proof it remains to note that any Sturmian sequence has
exactly one palindrome of length $k$ if $k$ is even \cite{DroPir}, and 
exactly one palindrome of length $k$ whose central letter is $0$ if $k$
is odd: the proof of Proposition 6 of \cite{DroPir} shows there is a
bijection between the set of palindromes of length $k+2$ and the set of
palindromes of length $k$ and that this bijection consists of erasing
the first and last (identical) letters. Since there is exactly one
palindrome of length $1$ with central letter $0$ and one palindrome of
length $1$ with central letter $1$, we see immediately that, for each odd
$k$, any Sturmian sequence has exactly one palindrome of length $k$ whose
central letter is $0$ and exactly one palindrome of length $k$ whose central
letter is $1$. \endpf

\begin{remark}\label{counterex}

\mbox{ }

\begin{itemize}

\item Another way of studying palindromes occurring in a sequence is to use 
the associated dynamical system and its (geometric) symmetry properties.
In this direction, the reader can look at \cite{AleBer} (see for example 
Theorem~19 in that paper) and \cite{BerVui} (where, in particular, a 
two-dimensional version of Theorem~\ref{roteseq} above is given).

\item Taking the first difference sequences modulo $2$ of $0,1$-sequences 
for studying factors and complexities of certain sequences was already 
used in \cite{Allfold,AllBou,Rot,Allschro} for example.

\item What is the palindrome complexity of {\sl any} sequence of complexity 
$2k$? Rote proves in \cite{Rot} that the fixed point $u$ of the morphism 
$0 \to 001$, $1 \to 111$ has complexity $2k$. With a slight modification of 
the argument in Theorem~\ref{general} below (this morphism is not primitive, 
but as soon as a word contains a $0$, we essentially know from which word it 
``comes by the morphism'') the reader can prove that ${\rm pal}_u(k) = 2$ 
for each $k \geq 1$. On the other hand, looking at another example of a
sequence of complexity $2k$ given in \cite{Rot}, namely the image under the 
morphism $a \to 0$, $b \to 1$, $c \to 10110$, $d \to 101$ of the fixed point 
of the morphism $a \to ad$, $b \to bac$, $c \to bacab$, $d \to baca$, it is 
easy to check that this sequence contains only one palindrome of length $6$.
(More precisely, we have for this sequence ${\rm pal}(k) = 2$, for $k = 1,
2, 3, 4, 5$, ${\rm pal}(k) = 1$, for $k = 6, 7, 8, 9, 10$, and 
${\rm pal}(k) = 0$ for $k \geq 11$.)

\item Let $v = v_0 v_1 \ldots $ be the fixed point of the morphism
$0 \to 001$, $1 \to 101$. The complexity of this sequence is given by
$\mbox{\rm fac}_v(k) = 2k$ for $k \geq 1$ (see~\cite{Fer1}). The reader can 
prove that the palindrome complexity of this sequence is $2$ for $k \leq 7$
and $0$ for $k \geq 8$. (Hint: prove there is no palindrome of length $8$ 
nor of length $9$, either by mimicking the method of Theorem~\ref{RS} above, 
or by using the recursive definition of the sequence $v$: $u_{3n} = u_n$, 
$u_{3n+1} = 0$, $u_{3n+2} = 1$.) In the same vein, the reader can prove that 
the Chacon sequence defined as the infinite fixed point of $0 \to 0010$, 
$1 \to 1$ does not contain palindromes of length $13$, nor palindromes of 
length $14$, hence it does not contain palindromes of length $\geq 13$. 
Note that the (usual) complexity of the Chacon sequence is equal to $2k-1$ 
for $k \geq 2$ (see~\cite{Fer1}).

\item We finally note that the (conjectured) palindrome complexity of
the Kolakoski sequence is also constant and equal to $2$ 
\cite[Section 4.1.3]{Lad}. 
Recall that the Kolakoski sequence is the sequence 
$$
2 \ 2 \ 1 \ 1 \ 2 \ 1 \ 2 \ 2 \ 1 \ \ldots
$$ 
defined as {\em the\,} sequence on the alphabet $\{1, 2\}$ that begins 
in $2$ and such that the sequence of its runlengths is equal to the sequence
itself (see \cite{Kol}, see \cite{Dek} for a recent survey).
Recall that a word $w$ on the alphabet $\{1,2\}$ is said to be differentiable 
if neither $111$ nor $222$ occurs, and its derivative $w'$ is the finite
sequence of lengths of blocks in $w$, discarding the first and/or last block 
if it has length one: $(12211)'=22$, $(121)'=1$. All factors of Kolakoski are 
{\em $C^{\infty}$-words}, i.e., words that can be differentiated infinitely
many times. The converse is an open conjecture. It is proved in 
\cite[Section 4.1.3]{Lad} that the number of $C^{\infty}$-palindromes of each 
length $\geq 1$ is $2$. Hence the palindrome complexity of the Kolakoski 
sequence is bounded by $2$, and it is conjecturally constant and equal 
to $2$.

\end{itemize}
\end{remark}

\section{Fixed points of uniform primitive morphisms}\label{genpal}

We first recall the definition of {\it class P morphisms} introduced by
Hof, Knill, and Simon in \cite{HKS}.

\begin{definition}[Hof-Knill-Simon]\label{defhks}
A morphism $\sigma$ on the (finite) alphabet ${\cal A}$ belongs to class P 
if there exists a palindrome $p$ and for every $a \in {\cal A}$ 
a palindrome $q_a$ such that, for every $a \in {\cal A}$, we have 
$\sigma(a) = p q_a$ (or, for every $a \in {\cal A}$, $\sigma(a) = q_a p$). 
The word $p$ can be empty. If $p$ is not empty, then some (or even all) 
$q_a$'s are allowed to be empty.
\end{definition}

We show in the following theorem that it is possible to compute the 
palindrome complexity of fixed points of uniform morphisms that belong
to class P.

\begin{theorem}\label{general}
Let $\sigma : {\cal A} \rightarrow {\cal A}^*$ be primitive and 
$u \in {\cal A}^{\Bbb N}$ such that $\sigma(u)=u$. We assume the following:

\begin{itemize}
\item[(i)] The morphism $\sigma$ belongs to class P: there exists a 
palindrome $p$ and, for every $a \in {\cal A}$, a palindrome $q_a$ such that  
$\sigma(a) = p q_a$ (the case where  $\sigma(a) = q_a p$ for every $a \in 
{\cal A}$ is analogous). 
\item[(ii)] The morphism $\sigma$ is uniform. Let $l_p := |p|$ and
$l_q := |q_a|$, for every $a \in {\cal A}$.
\item[(iii)] For $a \not= b$, $q_a$ and $q_b$ have distinct first (and hence 
last) symbols.
\end{itemize}  

Under these conditions we have the following recursion for the palindrome
complexity
$$
\exists n_0 \in {\Bbb N} \setminus \{0\}, \  \forall n \geq n_0, \ \ \
\mbox{\rm pal}(n) = \sum_{k \in {\cal E}} \mbox{\rm pal}(k),
$$
where
$$
{\cal E} = \{s, \ n = sl + l_p - 2j, \ \mbox{\rm where} \ 0 \leq j \leq l-1 \},
\ \mbox{\rm with } l:= l_p + l_q.
$$
\end{theorem}

\proof

Let us first explain our idea intuitively. Given any palindrome
$x = x_1 \ldots x_n$ occurring in $u$, we can associate with $x$ the
palindrome $\sigma(x) p = p q_{x_1} p q_{x_2} \ldots p q_{x_n} p$ which
is a factor of $u$ since $\sigma(u) = u$. Moreover, we can simultaneously
delete symbols at the beginning and end of $\sigma(x) p$, retaining the
palindromic form, without deleting  $q_{x_1}$ (and hence $q_{x_n}$)
completely. Thus we can associate with $x$ a family of palindromes.
By assumption (iii), this map is one-to-one when restricted to the set of
palindromes of length $n$ occurring in $u$. Conversely, we can do an inverse
procedure by re-substituting a given palindrome. That is, given a palindrome
we can decompose it (essentially uniquely by a result of Moss\'e
\cite{Mos1,Mos2}) and consider its inverse image under $\sigma$ which will
turn out to be a palindrome. Again by assumption (iii), this process is
one-to-one and we can thus establish a bijection between suitable sets of
palindromes. As a consequence, we get equality of their cardinalities
which yields recursive relations for the palindrome complexity of $u$.

Let us be more precise. Given $n \in {\Bbb N} \setminus\{0\}$ (this will
be the length of the long palindrome), we are looking for solutions $s$ in
${\Bbb N} \setminus \{0\}$ (this will be the length of the short palindrome
the long palindrome is coming from) of the equation
\begin{equation}\label{corresp}
n = s l + l_p - 2j \ \mbox{\rm subject to the condition } 0 \le j \le l - 1,
\ \mbox{\rm with } l = l_p + l_q.
\end{equation}
It is clear that there are at most two solutions to this problem, and if there
are two solutions $s_1,s_2$, then $l$ is even and $|s_1 - s_2| = 1$. We will
show the following:

{\it There exists $n_0 \in {\Bbb N} \setminus \{0\}$ such that for every 
$n \geq n_0$, we have}
\begin{equation}\label{rec}
{\rm pal}(n) \, = \sum_{{\tiny k \mbox{ solves Equation~\ref{corresp}}}}
{\rm pal}(k).
\end{equation}
Let ${\rm Pal}(m)$ denote the set of palindromes of length $m$ that are
factors of $u$. We therefore have ${\rm pal}(m) = \# {\rm Pal}(m)$.
Let $n \in {\Bbb N} \setminus \{0\}$. Assume first that Equation~\ref{corresp}
has two solutions $k,k+1$. We will define two maps
$$
\Phi : {\rm Pal}(k) \cup {\rm Pal}(k+1) \rightarrow {\rm Pal}(n), \;
\Psi : {\rm Pal}(n) \rightarrow {\rm Pal}(k) \cup {\rm Pal}(k+1)
$$
and show that, for $n$ large enough, they are one-to-one. Since ${\rm Pal}(k)$
and ${\rm Pal}(k+1)$ are clearly disjoint, this implies Equation~\ref{rec}.

The map $\Phi$ is defined as explained above. Start with a palindrome $x$
in ${\rm Pal}(k) \cup {\rm Pal}(k+1)$ and consider the associated palindrome
$\sigma(x) p$. Then, by Equation~\ref{corresp}, $|\sigma(x) p| \ge n$ and
we obtain an element of ${\rm Pal}(n)$ by ``pruning'' $\sigma(x) p$
suitably. Call this element $\Phi(x)$. It follows from assumption (iii) and
Equation~\ref{corresp} that $\Phi$ restricted to either ${\rm Pal}(k)$ or
${\rm Pal}(k+1)$ is one-to-one.
To see that $\Phi$ is one-to-one also on their union, we invoke the
recognizability property of $u$ as proven by Moss\'e, see \cite{Mos1,Mos2}.
This gives the claim immediately in the case $l_p \not= l_q$ (for $n$ large
enough to apply Moss\'e's result and to compare two decompositions) and
in the case $l_p = l_q$ one uses that $u$ is $N$-th power-free \cite{Mos1}
for some $N$ together with this argument to prove the claim.

The definition of $\Psi$ is slightly more complicated.
Let $x \in {\rm Pal}(n)$ be given. According to Moss\'e's result, for $n$
large enough, one can draw one ``bar'' and hence obtain a unique decomposition
of $x$ using assumption~(ii). (Recall that a sequence $u=(u_n)_{n \geq 0}$
fixed point of a morphism $\varphi$ can be written
$$
u = \varphi(u) = /\, \varphi(u_0) \, /\, \varphi(u_1) \, /\, \ldots
$$
and any factor of $u$ will ``contain'' bars -- possibly in a non-unique
way -- according to its position(s) in the sequence $u$. Note that for
a constant-length morphism, knowing the position of one bar gives the
positions of the other bars.) We will consider a slightly different 
decomposition. Namely, we will also draw bars between the $p$'s and the 
$q_a$'s, $a \in {\cal A}$. It is easy to see that for $n$ large enough, 
these modified decompositions are unique as well.

We thus have
\begin{equation}
x = \sigma \,/\, w_1 \,/\, w_2 \,/\, \ldots \,/\, w_r \,/\, \pi,
\end{equation}
where $w_i$ is equal to either $p$ or one of the $q_a$'s, $1\le i \le r$, 
$\sigma$ is a non-empty suffix and $\pi$ a non-empty prefix of $p$ or of 
some $q_a$. Since $x$ is a palindrome, the ``reflected modified decomposition''
\begin{equation}
x = \widetilde{\pi} \,/\, w_r \,/\, w_{r-1} \,/\, \ldots \,/\, w_1 \,/\, 
\widetilde{\sigma}
\end{equation}
corresponds to another modified decomposition of $x$ which must be the same by
uniqueness. This shows that $\sigma = \tilde{\pi}$ and $w_i = w_{r-i+1}$, 
$1 \le i \le r$. Moreover, the center of $x$ must also be a center of some 
$w_j$ in this representation. Finally, we have $w_i = p$ either for all 
odd $i$ or for all even $i$. Let us consider the first case. We have
\begin{equation}
x = \sigma \,/\, p \,/\, q_{x_1} \,/\, p \,/\, q_{x_2} \,/\, p \ldots \,/\, 
p \, /\, q_{x_2} \,/\, p \,/\, q_{x_1} \,/\, p \,/\, \pi
\end{equation}
and we can associate some $q_{x_0}$ with both $\sigma$ and $\pi$ 
(by assumption~(iii)). We therefore get the palindrome $x_0 x_1 x_2 \ldots x_2 
x_1 x_0$ which belongs to either $P(k)$ or $P(k+1)$, for otherwise 
Equation~\ref{corresp} would have a solution different from $k, k+1$. 
In the second case we proceed similarly, in this case $\sigma$/$\pi$ are 
suffix/prefix of $p$ and we obtain a palindrome $x_1 x_2 \ldots x_2 x_1$ 
which, by the same reasoning, belongs to $P(k)$ or $P(k+1)$. Call the obtained 
palindrome $\Psi(x)$. Using assumption~(iii) and the unique decomposition 
property, we see that $\Psi$ is one-to-one. This establishes the existence 
of a bijection between ${\rm Pal}(k) \cup {\rm Pal}(k+1)$ and ${\rm Pal}(n)$ 
and hence proves the assertion of the theorem in the case where 
Equation~\ref{corresp} has two solutions.

If there is one solution to Equation~\ref{corresp}, we can prove
Equation~\ref{rec} along the same lines, parts of the argument being even
simpler than in the two-solution case.

On the other hand, if ${\rm Pal}(n)$ is non-empty, we can define $\Psi$ as
above and our re-substitution argument then shows that Equation~\ref{corresp}
must have a solution. Hence the absence of such a solution implies
${\rm Pal}(n) = \emptyset$, that is, ${\rm pal}(n) = 0$. \endpf

\begin{remark}
This theorem applies to several prominent examples, such as the period
doubling morphism (${\cal A} = \{a,b\}$, $p = a$, $q_a = b$, $q_b = a$)
and the {\em square\,} of the Thue-Morse morphism (${\cal A} = \{a,b\}$,
$p = \varepsilon$, $q_a = abba$, $q_b = baab$). It also applies to
the square of some generalizations of the Thue-Morse morphism given
in \cite{AllSha2} (see Lemma 2 and Theorem 8 of \cite{AllSha2} in the case
where $m \in \{1, 2\}$).
\end{remark}

\bigskip

We know from Theorem~\ref{damzar} that the palindrome complexity
of a fixed point $u$ of a primitive morphism takes only finitely many
values. It then makes sense to ask whether the sequence 
$({\rm pal}_u(k))_{k \geq 1}$ itself is generated by a morphism.
If the morphism belongs to class P, satisfies some technical conditions,
and has constant length, we give an answer in our next theorem.

\begin{theorem}
Let $u = (u_n)_{n \geq 0}$ be a sequence that is a fixed point of a primitive
morphism belonging to class P and satisfying the extra conditions (i), (ii), 
(iii) of Theorem~\ref{general}. In particular the morphism is uniform: let $d$ 
be its length. Then the palindrome complexity $({\rm pal}_u(k))_{k \geq 1}$ is 
a $d$-automatic sequence.
\end{theorem}

\proof

Let $i \in [0, l-1]$ be a fixed integer and let $n$ be an integer $\geq 0$
($n \geq 1$ if $i=0$). We want to compute ${\rm pal}_u(ln+i)$ using 
Theorem~\ref{general} above. We distinguish two cases:

\begin{itemize}

\item If $l$ is odd, the congruence $l_p -2j \equiv i \bmod l$ has a unique
solution, say $j_0$, belonging to $[0, l-1]$. Then, from Theorem~\ref{general},
we have for $n$ large enough
$$
\mbox{\rm pal}_u(ln+i) = 
\mbox{\rm pal}_u \left(n + \frac{i-l_p+2j_0}{l}\right).
$$
Let $f(n):= {\rm pal}_u(n-1)$ for $n \geq 2$, then
\begin{equation}\label{kernel}
f(ln+i+1) = f\left(n + \frac{l+i-l_p+2j_0}{l}\right).
\end{equation}
We claim this implies that for $a \in [0,l-1]$, the sequence 
$(f(ln+a))_{n \geq \max(n_0,2)}$ is a linear combination of
the sequences $(f(n+p))_{n \geq \max(n_0,2)}$, with $p \in [0,3]$.

This is clear for $a \in [1,l-1]$ from Equation~\ref{kernel} above.
Using Equation~\ref{kernel} with $i = l-1$, we obtain for $n$ large enough,
$$
f(ln+l) = f\left(n + 1 + \frac{l-1-l_p+2j_0}{l}\right).
$$
Hence replacing $n$ by $n-1$,
$$
f(ln) = f\left(n + \frac{l-1-l_p+2j_0}{l}\right),
$$
and the claim is proved.

\item If $l$ is even, the congruence $l_p -2j \equiv i \bmod l$ either has no
solution (if $l_p$ and $i$ have opposite parities) or has two solutions, say 
$j_1$ and $j_2$, belonging to $[0, l-1]$. 
In the first case ${\rm pal}_u(ln+i)=0$. In the second case we have
for $n$ large enough
$$
\mbox{\rm pal}_u(ln+i) = 
\mbox{\rm pal}_u\left(n + \frac{i-l_p+2j_1}{l}\right)
+
\mbox{\rm pal}_u\left(n + \frac{i-l_p+2j_2}{l}\right).
$$
We conclude as above that $f(n):= {\rm pal}_u(n-1)$ for $n \geq 2$ has the 
property that for every $a \in [0,l-1]$, the sequence
$(f(ln+a))_{n \geq \max(n_0,2)}$ is a linear combination of
the sequences $(f(n+p))_{n \geq \max(n_0,2)}$, with $p \in [0,3]$.

\end{itemize}

Now from a result of \cite{AllSha3} quoted below as Theorem~\ref{allsha}, 
we see that the sequence $(f(n))_n$ (and hence the sequence 
$({\rm pal}_u(n))_{n \geq 1}$) is $d$-regular in the sense of \cite{AllSha1}. 
But this sequence takes only finitely many values, hence as noted in 
\cite{AllSha1} it must be $d$-automatic. \endpf

\bigskip

Before giving the next statement for the sake of completeness, we first 
recall (see \cite{AllSha1}) that a sequence $u = (u_n)_{n \geq 0}$ 
with values in a Noetherian ring $R$ is called $d$-regular for some 
integer $d \geq 1$ if the $R$-module generated by all subsequences
$(u_{d^t n + \ell})_{n \geq 0}$ with $t \geq 0$ and $\ell \leq d^t - 1$
has finite dimension. The following result is proved in \cite{AllSha3}.

\begin{theorem}[Allouche-Shallit]\label{allsha}
Let $u = (u_n)_{n \geq 0}$ be a sequence with values in a Noetherian ring $R$.
Suppose there exist an integer $d \geq 2$, an integer $t \geq 0$, an integer 
$r \geq 0$ and an integer $n_0 \geq 0$ such that each sequence
$(u_{d^{t+1} n + \ell})_{n \geq n_0}$ for $\ell \in [0, d^{t+1} - 1]$ is a 
linear combination of the sequences $(u_{d^j n + i})_{n \geq n_0}$ with 
$j \leq t$, $i \leq d^j -1$ and of the sequences $(u_{n + p})_{n \geq n_0}$ 
with $p \leq r$. Then the sequence $u$ is $d$-regular.
\end{theorem}

\section{Bounding the palindrome complexity in terms of the usual complexity}

Looking at some of the examples above we can ask several questions.

\begin{itemize}

\item We saw that Sturmian sequences, complement-symmetric Rote sequences, 
the (nonconstant) fixed point of $0 \to 001$, $1 \to 111$, and fixed points 
of primitive morphisms all have bounded palindrome complexity. On the other 
hand, for all these sequences the complexity satisfies $\mbox{\rm fac}(k) 
= O(k)$. Is it true that the property $\mbox{\rm fac}(k) = O(k)$ implies 
$\mbox{\rm pal}(k) = O(1)$? We will answer this question positively in 
Theorem~\ref{cassaigne} below. Note that the converse is not true, by far.
For example, apply to a binary sequence $u$ with complexity 
$\mbox{\rm fac}_u(k) = 2^k$ the morphism $0 \to 011001$, $1 \to 001011$. Then 
the complexity $\mbox{\rm fac}_v$ of the image $v$ of $u$ under this morphism
satisfies $2^{k/6} \leq \mbox{\rm fac}_v(k) \leq 9 \cdot 2^{k/6}$, whereas the 
palindrome complexity drops to $0$ from $k=8$.

\item More generally, for which sequences is it true that 
$\mbox{\rm pal}(k) = O(\mbox{\rm fac}(k)/k)$? For which sequences is it 
true that 
$0<\limsup_{k \to \infty}(k\mbox{ \rm pal}(k)/\mbox{\rm fac}(k)) < +\infty$?

\noindent
Of course this last assertion is not true for all sequences: take the sequence 
of binary digits of a normal real number. The complexity of this sequence is
$\mbox{\rm fac}(k) = 2^k$ and its palindrome complexity is
$\mbox{\rm pal}(k) = 2^{\lfloor (k+1)/2 \rfloor}$.
This is not true for the Rudin-Shapiro sequence (Theorem~\ref{RS})
nor for the sequence given in the third item of Remark~\ref{counterex}.
For this question see also Remark~\ref{remcor} below.

\item Is there a Pansiot-like theorem \cite{Pan} for the palindrome
complexity of fixed points of morphisms? Combining with the previous question,
is it true that, if the palindrome complexity of a fixed point of a morphism
is not ultimately $0$, then it satisfies
$$
0 < \limsup_{k \to \infty} (\mbox{\rm pal}(k)/\varphi(k)) < +\infty,
$$
where $\varphi(k)$ is either $1$, $\log \log k$, $\log k$ or $k$? 

This question can be first experimentally addressed by looking at all 
palindromes of reasonable lengths that occur in fixed points of various 
morphisms, getting an idea of what the order of magnitude of the 
palindrome complexity seems to be, and ... proving it. We give two
examples. 

\begin{itemize}

\item we know that the complexity of the (infinite) fixed point $u$ of the 
morphism $0 \to 001$, $1 \to 1$ satisfies $C_1 k^2 \leq \mbox{\rm fac}_u(k) 
\leq C_2 k^2$ \cite{Pan}. If we look for {\em maximal palindromes\,} in the 
sequence, i.e., factors of the sequence that are palindromes and that cannot 
be extended to longer palindromes occurring in the sequence, we find they 
are given by the recurrence $w_0 = 0$, $w_{m+1} = 1.\sigma(w_m)$. (Note that
by ``the palindrome $w$ can be extended to a palindrome occurring in the
sequence'' we mean of course that there exists a letter $a$ such that $awa$ 
occurs in the sequence.) It can be proved that $\mbox{\rm pal}(k)$ is the 
number of integers $m$ such that $m \leq k+1$, $|w_m| = 2^{m+1}+m-1 \geq k$, 
and $m-k$ is odd, which implies 
$\mbox{\rm pal}(k) = k/2 - 1/2 \log_2(k) + O(1)$.

\item In the same vein, the complexity of the fixed point beginning in $0$ of 
$0 \to 010$, $1 \to 11$ satisfies $\mbox{\rm fac}(k) \sim k \log_2\log_2(k)$ 
(see \cite{Cas2}). The behavior of the palindrome complexity of this sequence 
is different for odd and even length: $\mbox{\rm pal}(2n+1) = 1$ for 
$n \geq 4$, whereas $\mbox{\rm pal}(2n) \sim \log_2\log_2(2n)$.

\end{itemize}

\end{itemize}

Before we state the main theorem of this section, we need a definition and 
a preliminary result.

\begin{definition}

Let $w$ be a palindrome on the alphabet ${\cal A}$, and let $T$ be 
its period.

\begin{itemize}

\item If $T > |w|/2$, the palindrome $w$ is called {\em non-periodic}.

\item If $T \leq |w|/2$ and $T$ is odd, the palindrome $w$ is called
{\em a palindrome of odd period}.

\item If $T \leq |w|/2$ and $T$ is even, the palindrome $w$ is called
{\em a palindrome of even period}.

\end{itemize}

\end{definition}

\begin{remark}\label{perpal}

In this paper we thus call {\em periodic\,} a palindrome whose period (even
or odd) is at most half of the length of the palindrome. This means in 
particular that a periodic palindrome is a palindrome that can {\em a priori} 
be written as $A^d B$, where $B$ is a prefix of $A$ with $|B| < |A|$, and 
necessarily $d \geq 2$. It is not hard to see that the word $w$ is a
periodic palindrome if and only if there exist two palindromes $B$ and $C$, 
and an integer $d \geq 2$, such that $w = (BC)^d B$ (in the decomposition
$w = A^d B$, remember that $B$ is a prefix of $A$, put $A = BC$, and compare
the prefixes of length $|B|+|C|$ of $w$ and $\widetilde{w}$).

\end{remark}

\begin{lemma}\label{twin}

Let $w$ be a palindrome of even period $T$. The word $w$ is hence a
prefix of an infinite sequence $(xyxyxy \ldots)$ of smallest period $T$,
with $|x| = |y| = T/2 \leq |w|/4$. The prefix of length $|w|$ of the 
sequence $(yxyxyx \ldots)$ is denoted by $w^{\Theta}$ and called the 
{\em twin} of $w$. Then

\begin{itemize}

\item the word $w^{\Theta}$ is a periodic palindrome and its period is $T$. 
Furthermore $w^{\Theta} \neq w$;

\item ``taking the twin'' is an involution. More precisely the map 
$w \to w^{\Theta}$ is an involution on the set of palindromes of even
period. 

\end{itemize}

\end{lemma}

\proof

It is clear that $w$ and $w^{\Theta}$ have the same period, and that
the map $w \to w^{\Theta}$ is an involution. Note that $w \neq w^{\Theta}$, 
otherwise $x=y$ and $T/2$ would be a period of $w$. 

It remains to prove that $w^{\Theta}$ is a palindrome.
Let us write $w = (xy)^d z$, with $|x| = |y| = T/2$, $d \geq 2$, and
$z$ prefix of $xy$ with $|z| < |xy|$. Write $xy = zt$ (with $t$ non-empty).
As in Remark~\ref{perpal}, the words $z$ and $t$ must be palindromes.
Hence $t z = \widetilde{t} \widetilde{z} = \widetilde{zt} = \widetilde{xy}
= \widetilde{y} \widetilde{x}$. Hence the word $yz\widetilde{y}$, which is
a prefix of $yz\widetilde{y}\widetilde{x}$, is a prefix of $yztz = yxyz$,
hence a prefix of $yxyxy$. This shows that the word 
$y (xy)^{d-1} z \widetilde{y} = (yx)^{d-1}y z \widetilde{y}$ is a prefix 
of the sequence $(yxyxyx\ldots)$. Since the length of 
$y(xy)^{d-1} z \widetilde{y}$ is equal to $|w|$, we thus have an alternative 
definition of $w^{\Theta}$:

if $w$ is a palindrome of even period $T$ ($T \leq |w|/2$), 
with $w = (xy)^d z$, where $|x| = |y| = T/2$, $d \geq 2$, 
$z$ prefix of $xy$, and $|z|< |xy|$, then 
$$
w^{\Theta} = y (xy)^{d-1} z \widetilde{y}.
$$
Now, with the notations above, and remembering that $z$ and $t$ are
palindromes,
$$
\widetilde{w^{\Theta}} = y z (\widetilde{xy})^{d-1} \widetilde{y}
= y z (\widetilde{zt})^{d-1} \widetilde{y} =
y z (tz)^{d-1} \widetilde{y} =  y (zt)^{d-1} z \widetilde{y}
= y (xy)^{d-1} z \widetilde{y} = w^{\Theta}. \mbox{\endpf}
$$

\bigskip

We now give a theorem that bounds the palindrome complexity in terms of
the usual complexity, and that answers positively the first question
at the beginning of this section.

\begin{theorem}\label{cassaigne}
Let $u = u_0 u_1 u_2 \ldots$ be an infinite non-ultimately periodic sequence
on a finite alphabet. Then, for all $k \geq 1$, we have
$$
\mbox{\rm pal}_u(k) <
\frac{16}{k} \ \mbox{\rm fac}_u\left(k+\lfloor \frac{k}{4}\rfloor\right).
$$
\end{theorem}

\proof

We first suppose that the sequence $u$ is recurrent or that it is 
indexed by ${\Bbb Z}$ ($u = \ldots u_{-2} u_{-1} u_0 u_1 u_2 \ldots$).
In the latter case we suppose the sequence is not ultimately periodic 
``on the right''. Let $k \geq 1$ be a fixed integer. We split 
${\rm Pal}_u(k)$ the set of palindromes of length $k$ occurring in 
the sequence $u$ into three classes according to their periods $T$: 
$$
\begin{array}{ll}
{\rm Pal}_u^{(0)}(k) &:= \{w \in {\rm Pal}_u(k), \ T > k/2\}, \\

&\mbox{\rm (this is the set of non-periodic palindromes of length $k$
occurring in $u$)}, \\

\\

{\rm Pal}_u^{(1)}(k) &:= \{w \in {\rm Pal}_u(k), \ T \leq k/2 \ 
\mbox{\rm and $T$ odd}\}, \\

&\mbox{\rm (this is the set of palindromes of length $k$ and of odd period
occurring in $u$)}, \\

\\

{\rm Pal}_u^{(2)}(k) &:= \{w \in {\rm Pal}_u(k), \ T \leq k/2 \ 
\mbox{\rm and $T$ even}\}, \\

&\mbox{\rm (this is the set of palindromes of length $k$ and of even period
occurring in $u$)}. \\

\end{array}
$$

\begin{itemize}

\item For each $w \in \mbox{\rm Pal}_u^{(0)}(k)$, i.e., for each non-periodic 
palindrome $w$ of length $k$ in the sequence $u$ we choose an index $\ell$ 
(assuming furthermore, which is possible, that $\ell > \lfloor k/4 \rfloor$ 
if the sequence is recurrent) such that an occurrence of $w$ in $u$ begins 
at index $\ell$. Then, we associate with $w$ the language $S(w)$ that 
consists of the ($\lfloor k/4 \rfloor + 1$) words of length 
($k + \lfloor k/4 \rfloor$) that occur in $u$ and begin at indexes
between ($\ell - \lfloor k/4 \rfloor$) and $\ell$. (Note that these indexes
are nonnegative in the case where the sequence is recurrent, and well-defined 
if the sequence is indexed by ${\Bbb Z}$.) These words are all 
distinct; namely if the words beginning at indexes $\ell - i$ and $\ell - j$
(where $0 \leq i, j \leq \lfloor k/4 \rfloor$) were equal, the word $w$ would 
have a period $\leq |i - j| \leq k/4 \leq k/2$ (use Theorem~\ref{lynsch}), 
which contradicts the non-periodicity of $w$. We thus have
$$
\forall w \in {\rm Pal}_u^{(0)}(k), \ \ \ \ \ 
\#(S(w)) = \lfloor k/4 \rfloor + 1.
$$

\item For each $w \in \mbox{\rm Pal}_u^{(1)}(k)$, i.e., for each palindrome 
$w$ of length $k$ occurring in the sequence $u$ and having an odd period $T$, 
we choose an index $\ell$ (assuming furthermore, which is possible, that 
$\ell > \lfloor k/4 \rfloor$ if the sequence is recurrent) such that an 
occurrence of $w$ in $u$ begins at index $\ell$, and such that the word 
of length $|w|$ beginning at index $\ell + T$ is different from $w$ 
(this is possible since the sequence $u$ is not ultimately 
periodic). We define $S(w)$ as above. The words in $S(w)$ are pairwise 
distinct: if the words of $S(w)$ beginning at $\ell - i$ and $\ell - j$, 
say $z_i$ and $z_j$, (where $0 \leq i, j \leq \lfloor k/4 \rfloor$) were 
equal, then there would exist two words $A$ and $B$ of length $|j-i|$ such 
that $Az = zB$, where $z = z_i = z_j$. From Theorem~\ref{lynsch}
this implies the existence of two words $\gamma$ and $\delta$ and of an integer
$e \geq 0$ such that $A = \gamma \delta$, $B = \delta \gamma$, and 
$z = (\gamma \delta)^e \gamma$. Hence $Az = zB$ is a prefix of a periodic
sequence of period $|\gamma \delta| = |j-i|$. Since $w$ is a factor of
$Az = zB$ (this is the factor of length $k$ beginning at index $\ell$ in
the sequence $u$, and $Az = zB$ is the factor of length 
$k + \lfloor k/4 \rfloor$ beginning at index $\ell - \max(i,j)$ in the 
sequence $u$), this implies from Lemma~\ref{periods} that $Az$ has period
$T$ (note that $|j-i|+T \leq j + i + k/2 \leq 2 \lfloor k/4 \rfloor + k/2 
\leq k = |w|$). Now, the factor of the sequence $u$ of length $k$ beginning 
at index $\ell + T$ is a factor of $Az$ since $T \leq |j-i|$ (namely
$|j-i|$ is a period of $Az = zB$, which contains $w$, hence $|j-i|$ is 
a period of $w$, and by Lemma~\ref{periods} it is a multiple of $T$, so 
$T \leq |j-i|$) hence is equal to the factor of length $k$ beginning at 
$\ell$, i.e., to $w$, which gives a contradiction.
 
\item We consider now the palindromes of even period. For each 
$w \in \mbox{\rm Pal}_u^{(2)}(k)$, i.e., for each palindrome $w$ of length 
$k$ and of even period $T$ (remember that means in particular that 
$T \leq k/2$), we consider its twin $w^{\Theta}$ (see Lemma~\ref{twin}).
Since the map $w \to w^{\Theta}$ is involutive on the set of palindromes 
of even period, and $w^{\Theta} \neq w$ (Lemma~\ref{twin}), palindromes 
of even period can be grouped in pairs of twins. For each pair of twins 
$(w, w^{\Theta})$ of length $k$ and (even) period $T$, such that at least 
one of $w$, $w^{\Theta}$ occurs in the sequence $u$, we choose an index 
$\ell$ (and we can suppose $\ell \geq \lfloor k/4 \rfloor$ in the case where 
$u$ is recurrent) such that one of the words $w$, $w^{\Theta}$ occurs in $u$ 
at index $\ell$, and such that its twin does not occur at index $\ell + T/2$, 
this is possible since the sequence $u$ is not ultimately periodic: if it were 
true that for every index $\ell$ such that the word $w$ (resp.\ $w^{\Theta}$) 
occurs in the sequence $u$ at index $\ell$, then the word $w^{\Theta}$ 
(resp.\ $w$) would occur in the sequence $u$ at index $\ell + T/2$, then both 
words $w$ and $w^{\Theta}$ would occur in $u$, and, for every index $\ell$ 
such that $w$ occurs at $\ell$, $w$ would also occur at index $\ell + T$, 
which would imply that $u$ is ultimately periodic. We construct the set 
$S(w)$ (or $S(w^\Theta)$) as above, and the words in $S(w)$ (or $S(w^\Theta)$)
are again pairwise distinct.

\end{itemize}

We thus constructed at least $\#\mbox{\rm Pal}_u^{(0)}(k) +
\#\mbox{\rm Pal}_u^{(1)}(k) + \#\mbox{\rm Pal}_u^{(2)}(k)/2 
\geq \mbox{\rm pal}_u(k)/2$ languages, each of which contains 
$\lfloor k/4 \rfloor + 1$ words of length $k + \lfloor k/4 \rfloor$
that are distinct. We prove now that these sets are pairwise disjoint.
If there exists a word $z$ belonging to $S(w_1) \cap S(w_2)$ (hence 
$|z| = k + \lfloor k/4 \rfloor$), then there exist indexes $\ell_1$ and 
$\ell_2$, nonnegative integers $i_1$ and $i_2$ both $\leq \lfloor k/4 \rfloor$,
such that $w_1$ occurs at index $\ell_1$ in the sequence $u$, $w_2$ 
occurs at index $\ell_2$ in $u$, and $z$ occurs in $u$ both at indexes 
$\ell_1 - i_1$ and $\ell_2 - i_2$. Then $w_1$ occurs at position $i_1$ 
in $z$ and $w_2$ at position $i_2$ in $z$. Let $h := i_2 - i_1$ 
(assuming $i_2 \geq i_1$). If $h = 0$, then $w_1=w_2$.
Otherwise we have $h > 0$. Then the prefix of length $k - h$ of $w_2$ is
equal to the suffix of length $k - h$ of $w_1$. We write $w_1 = At$ and 
$w_2 = tB$, where $|t| = k-h$ and $|A|=|B|=h$. Now let $z' := AtB$ be the
factor of $z$ of length $k+h$ obtained by ``superposing'' $w_1$ and $w_2$.
Since $w_1$ and $w_2$ are both palindromes, we have:
$$
(A \widetilde{B}) \widetilde{t} = A (\widetilde{B} \widetilde{t} ) =
A(tB) =  z' = (At)B = (\widetilde{t} \widetilde{A}) B = 
\widetilde{t} (\widetilde{A} B).
$$
Hence, from Theorem~\ref{lynsch} there exist two words $C$ and $D$ and 
an integer $q \geq 0$ such that $\widetilde{t} = (CD)^q C$, $A \widetilde{B}
= CD$, and $\widetilde{A} B = DC$. Hence $z' = (CD)^{q+1}C$, and
$|CD| = |A \widetilde{B}| = 2h$ is a period of $z'$.  In particular $2h$ is 
a period of $w_1$. Since $2h \leq 2i_2 \leq 2 \lfloor k/4 \rfloor \leq k/2$,
the word $w_1$ is a periodic palindrome whose period, say $T$, divides $2h$:
this is a consequence of Lemma~\ref{periods} which also gives that $T$ is 
a period of $z'$ and hence of $w_2$. 

\begin{itemize}

\item If $T$ divides $h$, then the equality
$z' = A \widetilde{B} \widetilde{t}$ shows that $A = \widetilde{B}$
(remember that $|A| = |B| = h$). Hence $w_1 = At$ and 
$w_2 = tB = \widetilde{B} \widetilde{t}$ coincide on their
prefixes of length $|A| = h$, hence $w_1 = w_2$. 

\item If $T$ does not divide $h$, 
let $A \widetilde{B} = E^r$, with $|E| = T$ and $r \geq 1$. Since 
$2h = |A \widetilde{B}| = |E^r| = rT$ and $T$ does not divide $h$, $T$ must 
be even, and $r$ must be odd. Let $E = xy$, with $|x| = |y| = T/2$. 
We have $A \widetilde{B} = (xy)^r$, which implies, since $r$ is odd,
that $A = (xy)^{(r-1)/2} x$ and $\widetilde{B} = y (xy)^{(r-1)/2}$.
The word $xy$ is a prefix of $A \widetilde{B}$, hence of $A \widetilde{B}
\widetilde{t} = z'$. Since $w_1$ is a prefix of $z'$, and since 
$|w_1| \geq 2T = 2|xy|$, we see that $xy$ is a prefix of $w_1$. Finally 
$w_1$ is the prefix of length $k$ of the sequence $(xyxyxy\ldots)$ and the 
smallest period of this sequence is $T$.

Now $\widetilde{B} \widetilde{t}D = \widetilde{B}(CD)^{q+1} =
y (xy)^{r(q+1) + (r-1)/2}$. Hence $yx$ is a prefix of 
$\widetilde{B}\widetilde{t}D$. Since $w_2 = \widetilde{B}\widetilde{t}$ 
is a prefix of $\widetilde{B}\widetilde{t}D$, and $|w_2| \geq 2 |xy|$, the
word $yx$ is a prefix of $w_2$. Now $T = |xy|$ is a period of $w_2$, hence
we see that $w_2$ is the prefix of length $k$ of the sequence $(yxyxyx\ldots)$. 
Since the smallest period of this sequence must be $T$ (as for the sequence
$(xyxyxy\ldots)$), $w_2$ is the twin of $w_1$, which is a contradiction, 
since we chose exactly one representative in each pair of twins.

\end{itemize}

\noindent We thus have
$$
\frac{\mbox{\rm pal}_u(k)}{2} (\lfloor k/4 \rfloor + 1) \leq
\mbox{\rm fac}_u(k + \lfloor k/4 \rfloor).
$$
Since $\lfloor k/4 \rfloor + 1 > k/4$, we have that
$$
\mbox{\rm pal}_u(k) <
\frac{8}{k} \ \mbox{\rm fac}_u\left(k+\lfloor \frac{k}{4}\rfloor\right).
$$

\bigskip

Now, if we take a non-ultimately periodic sequence $u = (u_n)_{n \geq 0}$ on 
an alphabet ${\cal A}$ and if $u$ is not recurrent, let $\omega$ be a letter 
that does not belong to ${\cal A}$. Define the sequence 
$u^*=(u^*_n)_{n \in {\Bbb Z}}$ on ${\Bbb Z}$ by $u^*_n = u_n$ if $n \geq 0$ 
and $u^*_n = \omega$ for $n < 0$. We clearly have:
$$
\begin{array}{ll}
\mbox{\rm pal}_u(k) 
= 
\mbox{\rm pal}_{u^*}(k) - 1
&<
\displaystyle
\frac{8}{k} \ \mbox{\rm fac}_{u^*} (k+\lfloor\frac{k}{4}\rfloor)-1
=
\frac{8}{k} \ \left(\mbox{\rm fac}_u (k+\lfloor\frac{k}{4}\rfloor)
+ (k+\lfloor\frac{k}{4}\rfloor)\right) - 1
\\
&<
\displaystyle
\frac{16}{k} \ \mbox{\rm fac}_{u} 
\left(k+\lfloor\frac{k}{4}\rfloor\right). \mbox{\endpf}
\end{array}
$$

\begin{remark}\label{remcor}

\mbox{ }

\begin{itemize}

\item As a corollary of Theorem~\ref{cassaigne} above, we see that {\em a 
sequence such that $\mbox{\rm fac}(k) = O(k)$ has bounded palindrome 
complexity}. This result can be compared to a result of Cassaigne \cite{Cas1}:
{\em if the complexity of a sequence satisfies $\mbox{\rm fac}(k) = O(k)$,
then $(\mbox{\rm fac}(k+1)- \mbox{\rm fac}(k)) = O(1)$}.

\item In particular {\em any} automatic sequence has bounded palindrome 
complexity, and any fixed point of a primitive morphism has bounded complexity 
(thus recovering Theorem~\ref{damzar}). Namely the (block)-complexity of an 
automatic sequence is $O(k)$ \cite{Cob}, and the block complexity of a fixed 
point of a primitive morphism is also $O(k)$ \cite{Mic1,Mic2}.

\item In view of Theorem~\ref{cassaigne}, a natural question is: is there 
a universal upper bound for the quantity 
$k\mbox{ \rm pal}(k)/\mbox{\rm fac}(k)$? 
(The quantity $\mbox{\rm fac}(k + \lfloor k/4\rfloor)$ instead of 
$\mbox{\rm fac}(k)$ seems to appear only for technical reasons).
The answer is {\bf no}, as shown by the following example, for which
$k\mbox{ \rm pal}(k) / \mbox{\rm fac}(k)$ reaches $\sqrt k/4$ for 
certain values of $k$, while $\mbox{\rm fac}(k) = O(k^{3/2})$;
We only outline the proof. Start with $w_0 = 1$, and define 
$w_{j+1} := w_j x_1 x_2 \ldots x_{2^{2^j-1}}$
where  $x_i = 0^{2^{2^{j+1}}+2-4i} \widetilde w_j 0^{2^{2^{j+1}}-4i} w_j$.
Then define the sequence $u$ as the limit of $w_j$ when $j$ tends to 
infinity. We have
$$
w_0 = 1, \ \ w_1 = 10011, 
$$
$$
w_2 = 100110000000000000011001000000000000100110000000000110010000000010011
$$
$\ldots  \ \ (|w_3| = 4997)$.
It can be proved that
$$
|w_j| \sim 2^{3.2^{j-1}},  \ \ \
\mbox{\rm fac}_u(2^{2^j}) \sim 2^{2^j+2},  \ \ \
\mbox{\rm pal}_u(2^{2^j}) = 2^{2^{j-1}}+1.
$$ 
Hence $\log \mbox{\rm fac}_u(k) / \log k$ oscillates between $1$ and $3/2$, 
while $\log \mbox{\rm pal}_u(k) / \log k$ oscillates between $0$ and $1/2$, 
but these two quantities are not ``in phase'' so that their difference 
oscillates between $1/2$ and $3/2$. Note that it might be the case that
$\sqrt{k}\mbox{ \rm pal}(k) / \mbox{\rm fac}(k)$ is still universally bounded.

On the other hand, the constant $1/4$ in the theorem can be changed, so
that the quantity 
$k \mbox{ \rm pal}(k) / \mbox{\rm fac}(k + \lfloor \alpha k \rfloor)$ is 
universally bounded for any fixed $\alpha > 0$, the bound depending of 
course of $\alpha$.

\item C. Choffrut \cite{Cho} asked the following nice question: is it true
that 
$$
\mbox{\rm pal}_u(k) = O(\sqrt{\mbox{\rm fac}_u(k)})?
$$ 
Note that this is true for sequences $u$ with maximal complexity, or even for 
sequences with complexity $\Theta(\alpha^k k^{-a})$ for some $\alpha > 1$ and 
$a \geq 0$ (use the easy bound 
$\mbox{\rm pal}_u(k) \leq \mbox{\rm fac}_u(\lfloor \frac{k+1}{2} \rfloor)$. 
This is also true from Theorem~\ref{cassaigne} above if 
$\mbox{\rm fac}_u(k) = O(k^{3/2})$, or if $\mbox{\rm fac}_u(k) = \Theta(k^2)$.
Hence, using Pansiot's theorem \cite{Pan}, this is true for fixed points of 
non-trivial morphisms.
 
\end{itemize}

\end{remark}

\section{On a question of Hof, Knill, and Simon}

To end this paper, we recall a question of Hof, Knill, and Simon in 
\cite[Remark 3, p.\ 153]{HKS}: {\sl are there (minimal) sequences containing
arbitrarily long palindromes that arise from substitutions none of which
belongs to class P?} This question is still open. But we prove in 
Theorem~\ref{particularhks} below that we can restrict ourselves to particular 
substitutions in class P, as well as to nonperiodic (minimal) sequences.
We first prove the following result.

\begin{lemma}\label{lemmaparticularhks}

Let $u$ be a fixed point of a primitive morphism $\sigma$ on the alphabet
${\cal A}$. Suppose that there exists a non-empty word $x$ that is a prefix 
of $\sigma(a)$ for all $a \in {\cal A}$ (resp.\ a suffix of $\sigma(a)$ 
for all $a \in {\cal A}$). Write $\sigma(a) = x z_a$ 
(resp.\ $\sigma(a) = z_a x$). Let $\sigma_{\#}$ be the morphism defined 
by $\sigma_{\#}(a) = z_a x$ (resp.\ $\sigma_{\#}(a) = x z_a$). Then 
$\sigma_{\#}$ is primitive, and any fixed point $v$ of a power of 
$\sigma_{\#}$ (there always exists at least one such fixed point) has the 
same factors as $u$.

\end{lemma}

\proof

The morphism $\sigma_{\#}$ is clearly primitive. Suppose now, 
for example, that $\sigma(a) = x z_a$ for all $a \in {\cal A}$. 
Then, for all $a \in {\cal A}$ we have $x \sigma_{\#}(a) = \sigma(a) x$.
Hence, for any word $w$ on ${\cal A}$, we have 
$x \sigma_{\#}(w) = \sigma(w) x$. Taking for $w$ prefixes of increasing 
length of $u$, we get a limit when this length goes to infinity:
$x \sigma_{\#}(u) = \sigma(u) = u$. 
Hence, $x \sigma_{\#}(x) \sigma_{\#}^2(x) \ldots  = u$. From this equality 
we see that, calling $\alpha$ the first letter of $x$ (remember that $x$
is not empty), any factor of $u$ is a factor of $\sigma_{\#}^j(\alpha)$ for 
$j \geq j_0$ (remember that $u$ is minimal, since $\sigma$ is primitive). 
Trivially any factor of $\sigma_{\#}^j(\alpha)$ for any $j$ is a factor of 
$u$. Now there exist integers $\ell > 0$ such that $\sigma_{\#}^{\ell}$ 
admits a fixed point, say $v$. For any such $\ell$ it is clear from what 
precedes that $v$ and $u$ have the same factors. 
The case where $x$ is a suffix of $\sigma(a)$ for all $a$ is similar. 
\endpf

\begin{theorem}\label{particularhks}

\mbox{ }

\begin{itemize}

\item Let $u$ be a sequence that is a fixed point of a primitive morphism in 
class P. Then the set of factors of $u$ is the same as the set of factors
of a sequence $u_{\#}$ that is a fixed point of a primitive morphism in 
class P where the palindrome $p$ occurring in Definition~\ref{defhks}
has length $0$ or $1$. We can even choose one of the two forms
$a \to p q_a$ for all $a$, or $a \to q_a p$ fo all $a$.

\item Let $u$ be a {\em periodic} sequence that contains arbitrarily long 
palindromes, then $u$ is a fixed point of a morphism in class P.

\end{itemize}

\end{theorem}

\proof

Let $u$ be a sequence that is a fixed point of the primitive morphism $\sigma$ 
on the alphabet ${\cal A}$. Suppose there exists a palindrome $p$ and, for 
any $a \in {\cal A}$ a palindrome $q_a$, with $\sigma(a) = p q_a$.

\medskip

If $p$ is empty, the first assertion of Theorem~\ref{particularhks} is 
satisfied.

\medskip

If $p$ has even length, let $p = r \widetilde{r}$ for some word $r$. 
Define the morphism $\sigma_{\#}$ by $\sigma_{\#}(a) := \widetilde{r} q_a r$
for all $a \in {\cal A}$. Applying Lemma~\ref{lemmaparticularhks} with
$x=r$, we know that this morphism is primitive; furthermore there exists an 
integer $\ell > 0$ such that $\sigma_{\#}^{\ell}$ admits a fixed point, 
say $v$, and the sequences $v$ and $u$ have the same factors. Finally, the 
definition of $\sigma_{\#}$ shows that the image of any letter is a palindrome, 
hence that the image of any palindrome is also a palindrome. This proves that 
the image of any letter by $\sigma_{\#}^{\ell}$ is a palindrome. In other 
words, $\sigma_{\#}^{\ell}$ belongs to class P, and the corresponding 
palindrome $p$ is empty.

\medskip

If $p$ has odd length, write $p = r b \widetilde{r}$ for some letter $b$ and
some word $r$. Define, for all $a \in {\cal A}$, $\sigma_{\#}$ by 
$\sigma_{\#}(a) :=  b \widetilde{r} q_a r$. 
Applying Lemma~\ref{lemmaparticularhks} with $x = r$, we mimic the proof 
just above in the case where the length of $p$ is even, leading to
a morphism in class P, whose corresponding palindrome $p$ has length $1$.

\medskip

Finally, if the morphism has, for example, the form $a \to p q_a$ for all $a$, 
we can apply Lemma~\ref{lemmaparticularhks} with $x = p$ to obtain a morphism 
of the form $a \to q_a p$ for all $a$. 

\bigskip

We now prove that the answer to the question of Hof, Knill and Simon is 
negative for a periodic sequence (that is a fixed point of a primitive
morphism). Let $u = www\ldots$ be a periodic sequence that contains 
arbitrarily long palindromes. Let $s$ be a palindromic factor of $u$ such that 
$|s| \geq 2 |w|$. We can write $s = x w^k y$, where $k \geq 1$, and 
$0 \leq |x|, |y| < |w|$. Since $s$ is a palindrome, we have 
$s = \widetilde{y} (\widetilde{w})^k \widetilde{x}$. Hence $\widetilde{w}$ is
a factor of $s$, hence a factor of $u$. Thus $\widetilde{w}$ must be a factor
of $ww$. Let $ww = A \widetilde{w} B$. Since $|\widetilde{w}| = |w|$, we see
that $|w| = |A| + |B|$, then we must have $w = AB$ ($A$ is a prefix of
$w$ and $B$ is a suffix of $w$). The equality $ww = A \widetilde{w} B$ then
implies that $\widetilde{w} = BA$. Hence $w = \widetilde{A} \widetilde{B}$, 
which shows that $A$ and $B$ are palindromes. We conclude by noting that the 
sequence $u = w w \ldots$ is a fixed point of the morphism $\tau$ defined 
by, for all $a \in {\cal A}$,  $\tau(a) := w = AB$ that is in class P. \endpf

\bigskip

\noindent
{\bf Acknowledgments.} JPA, MB, and DD began to discuss together on 
palindromes during the summer school and workshop ``The Mathematics of 
Aperiodic Order'' at the University of Alberta at Edmonton. They want 
to thank very warmly R. V. Moody for having organized this conference
and having permitted many people from different fields to confront
questions and knowledges.
The second author thanks the German Research Council (DFG) for its
support. The fourth author thanks the German Academic Exchange 
Service for financial support through Hochschulsonderprogramm III 
(Postdoktoranden).
Finally we would like to thank A. Rodenhausen for her questions on the
Kolakoski sequence, S. Brlek for his helpful comments on this sequence, 
and V. Berth\'e for several useful comments on a previous version of this 
paper.

\end{document}